\documentclass[10pt,amssymb]{article}
\usepackage{amsfonts,amsmath,latexsym}
\title{Poincar\'e-Einstein metrics and \\ 
the Schouten  tensor}
\author{Rafe Mazzeo\thanks{Supported by the NSF under Grant 
DMS-9971975 and at MSRI by NSF grant DMS-9701755} \\ 
Stanford University \and  Frank Pacard\thanks{Supported at 
MSRI by NSF grant DMS-9701755} \\ Universit\'e Paris XII}
\date{}

\begin{document}
\maketitle

\newcommand{\CC}{\mathbb C}
\newcommand{\HH}{\mathbb H}
\newcommand{\RR}{\mathbb R}
\newcommand{\del}{\partial}
\newcommand{\e}{\epsilon}
\newcommand{\olg}{\overline{g}}
\newcommand{\calC}{{\mathcal C}}
\newcommand{\calE}{{\mathcal E}}
\newcommand{\calF}{{\mathcal F}}
\newcommand{\calG}{{\mathcal G}}
\newcommand{\calL}{{\mathcal L}}
\newcommand{\calM}{{\mathcal M}}
\newcommand{\calP}{{\mathcal P}}
\newcommand{\calS}{{\mathcal S}}
\newcommand{\calU}{{\mathcal U}}
\newcommand{\calV}{{\mathcal V}}
\newcommand{\Ric}{{\mathrm{Ric}}}
\newcommand{\tr}{{\mathrm{tr}}\,}
\newcommand{\Mbar}{\overline{M}}
\newcommand{\dMbar}{\del\Mbar}
\newcommand{\ska}{\sigma_k(A)}
\newcommand{\sk}{\sigma_k}

\newtheorem{theorem}{Theorem}
\newtheorem{proposition}{Proposition}
\newtheorem{corollary}{Corollary}
\newtheorem{lemma}{Lemma}
\newtheorem{definition}{Definition}
\newtheorem{remark}{Remark}

\begin{abstract}
We examine here the space of conformally compact metrics $g$ on the 
interior of a compact manifold with boundary which have 
the property that the $k^{\mbox{th}}$ elementary symmetric 
function of the Schouten tensor $A_g$ is constant. When $k=1$
this is equivalent to the familiar Yamabe problem, and the
corresponding metrics are complete with constant negative
scalar curvature. We show for every $k$ that the deformation
theory for this problem is unobstructed, so in particular
the set of conformal classes containing a solution of any one
of these equations is open in the space of all conformal classes.
We then observe that the common intersection of these solution spaces 
coincides with the space of conformally compact
Einstein metrics, and hence this space is a finite intersection
of closed analytic submanifolds.
\end{abstract}

Let $\Mbar^{\, n+1}$ be a smooth compact manifold with boundary. A metric 
$g$ defined on its interior is said to be {\it conformally compact} if 
there is a nonnegative defining function $\rho$ for $\dMbar$ (i.e.\  
$\rho = 0$ only on $\dMbar$ and $\nabla \rho \neq 0$ there)
such that $\olg = \rho^2 g$ is a nondegenerate metric on $\Mbar$.
The precise regularity of $\rho$ and $\olg$ is somewhat peripheral
and shall be discussed later. Such a metric is automatically complete. 
Metrics which are conformally compact and also Einstein are of great 
current interest in (some parts of) the physics community, since they 
serve as the basis of the AdS/CFT correspondence \cite{Wi}, and they 
are also quite interesting as geometric objects. Since they are natural 
generalizations of the hyperbolic metric on the ball $B^{n+1}$, as 
well as the complete constant negative Gauss curvature metrics 
on hyperbolic Riemann surfaces -- which exist in particular on the 
interiors of arbitrary smooth surfaces with boundary -- and which 
are often called Poincar\'e metrics \cite{MT}, we say that a metric 
which is both conformally compact and Einstein is {\it Poincar\'e-Einstein} 
(or P-E for short). Until recently, beyond a handful of examples, the only 
general existence result concerning the existence of P-E metrics was the 
local perturbation theory of Graham and Lee \cite{GL}, which gives 
an infinite dimensional family of such metrics in a neighbourhood of 
the hyperbolic metric on the ball, parametrized by conformal classes
on the boundary sphere near to the standard one. Recently many new
existence results have been obtained, including further perturbation
results by Biquard \cite{Bi} and Lee \cite{L}, 
and Anderson has announced some important global existence results in 
dimension four \cite{A3}. Many interesting geometric and topological 
properties of 
these metrics have also been found \cite{G}, \cite{XW}, \cite{A1}, \cite{A2};
this last paper also surveys a number of intriguing examples of P-E metrics. 

A common thread through the analytic approaches to the construction of these 
metrics is the possible existence of an $L^2$ obstruction, or more simply a 
finite dimensional cokernel of the (suitably gauged) linearization of the 
Einstein equations around a solution. For any P-E metric where this 
obstruction is trivial, the implicit function theorem readily implies that 
the moduli space $\mathcal E$ of P-E metrics is (locally) a Banach manifold, 
parametrized by conformal classes of metrics on $\dMbar$. (Actually, the 
smoothness of $\mathcal E$ is true in generality \cite{A3}, 
but this geometric parametrization breaks down.) Unfortunately, the only 
known geometric criteria ensuring the vanishing of this obstruction are 
strong global ones \cite{L}.

The purpose of this note is to introduce some new ideas into this 
picture which may help elucidate the structure of this moduli space.
We consider a related family of conformally compact metrics which
satisfy a (finite) set of {\it scalar} nonlinear equations, including and 
generalizing the familiar Yamabe equation, which we introduce below. 
These are sometimes called the $\sk$-Yamabe equations, $k=1, \cdots, n+1$.
The hyperbolic metric on the ball, or indeed an arbitrary P-E metric on 
any manifold with boundary satisfies each of these equations, and 
conversely, in this particular (conformally compact) setting, metrics 
which satisfy every one of these scalar problems are also P-E. The punchline 
is that the deformation theory for the $\sk$-Yamabe equations is always 
unobstructed! This fact seems to have been unappreciated even
in the setting of compact manifolds without boundary (except for
the case $k=1$). The full implications of this statement in the
conformally compact case for the moduli space of P-E metrics is 
not completely evident at this point, but this relationship seems
quite likely to be of some value. Furthermore, the deformation theory
for these $\sk$-Yamabe metrics is new, and also of some interest.

To define these equations, recall the Schouten tensor $A_g$, defined for 
any metric $g$ on a manifold of dimension $n+1$ by the formula
\begin{equation}
A_g = \frac{1}{n-1}\left(\Ric - \frac{R}{2n}g\right);
\label{eq:sch}
\end{equation}
here $\Ric = \Ric_g$ and $R = R_g$ are the Ricci tensor and scalar curvature 
function for $g$. This tensor occupies a prominent position in conformal 
geometry because it transforms quite nicely under conformal changes of 
metric. In fact, if $\tilde{g} = e^{2u}g$, then 
\begin{equation}
A_{\tilde{g}} = A_g - \nabla^2 u + du \otimes du - \frac12 |\nabla u|^2 g.
\label{eq:trsch}
\end{equation}
Notice also that $g$ is Einstein if and only if $A_g$ is diagonal:
\begin{equation}
\Ric_g = \frac{R}{n+1}g \Longleftrightarrow A_g = \frac{R}{2n(n+1)}g.
\label{eq:AEg}
\end{equation}

The metric $g$ is a $\sk$-Yamabe metric if the $k^{{\mathrm{th}}}$
elementary symmetric function of the eigenvalues of $A_g$ are constant,
$\sk(A_g) = C$. This is usually posed as a problem in conformal
geometry: starting with an arbitrary metric $g$, the $\sk$-Yamabe
problem is to find a new metric $\tilde{g} = e^{2u}g$ in the
conformal class of $g$ such that $\sk(A_{\tilde{g}}) = C$. 
Notice that $\sigma_1(A_g) = R/2n$, and so $\tilde{g}$ is a 
$\sigma_1$-Yamabe metric if and only if its scalar curvature
is constant. When $M$ is compact (without boundary), this equation
is semilinear (for the function $v$ defined by $v^{4/(n-1)} = e^{2u}$)
and the existence theory is complete and by now well known \cite{LP}.
However, when $k > 1$ these equations are fully nonlinear and the 
existence theory is much less well understood. Recent significant
progress has been made by Chang-Gursky-Yang \cite{CGY} when $k=2$,
and also by Viaclovsky \cite{V}, but much remains to be understood. 
In particular, in contrast with the ordinary 
Yamabe problem, for $k>1$ the $\sigma_k$-Yamabe problem seems
to be somewhat more tractable for positively curved metrics:
a crucial a priori ${\mathcal C}^2$ estimate is missing
in the negative case \cite{V}.   

We now write out the $\sigma_k$-Yamabe equations (within a conformal
class) more explicitly. Fixing $g$ and using (\ref{eq:trsch}),
we see that 
\begin{equation}
\tilde{g} = e^{2u}g \qquad \mbox{satisfies} 
\qquad \sk(A_{\tilde{g}}) = (-1)^k C
\label{eq:skyep}
\end{equation}
provided
\begin{equation}
\calF_k(u,C) = \sk(\nabla^2 u - du \otimes du + \frac12|\nabla u|^2 g
- A_g) - C e^{2ku} = 0. 
\label{eq:skye}
\end{equation}
The symmetric function of the eigenvalues of $A_{\tilde{g}}$ here
is computed with respect to $g$ rather than $\tilde{g}$, which accounts 
for the exponential factor; the sign on the final term comes from taking 
$\sk$ of $-A_{\tilde{g}}$. We have (momentarily) suppressed the dependence 
of the functional $\calF_k$ on $g$, although this will be important later. 
In fact we shall also regard $\calF_k$ as a functional on the
space of metrics as follows. Let $\mathfrak S$ be a (local) slice of the
set of conformal classes of $M$; in other words, we assume that
$\mathfrak S$ contains one metric from each conformal class in a small
neighbourhood in the space of all conformal classes. It is customary
when defining these slices also to mod out by the diffeomorphism group,
but we do not do this here as there is no particular need.
Then we may also consider 
$(g,u) \to \calF_k(g,u,C)$, for $g \in \calS$ and $u$ a (smooth) 
function on $M$. For any constant $\beta_k$, define
\begin{equation}
\Sigma_k(\beta_k) = \{(g,u): \calF_k(g,u,\beta_k) = 0\}.
\label{eq:defSigk}
\end{equation}
Clearly this set does not depend on the choice of slice $\mathfrak S$,
and so defines a subset within the space of all metrics on $M$. 

As already indicated, the main result here involves 
the perturbation theory for solutions of $\calF_k$, or equivalently,
the structure of the sets $\Sigma_k(\beta_k)$, in the 
case where $M^{n+1}$ is a manifold with boundary and all metrics
are conformally compact. We use the same notation, namely
$\Sigma_k(\beta_k)$, to denote the set of conformally compact
$\sk$-Yamabe metrics on $M$. Since any conformally compact metric 
has asymptotically negative (in fact, isotropic) sectional
curvatures, we see that of necessity $\beta_k > 0$. 
In particular, the particular constants $\beta_k$ corresponding
to the hyperbolic metric $g_0$ on $B^{n+1}$ are
\begin{equation}
\beta_k^0 = 2^{-k} { n+1 \choose k }.
\label{eq:bkhyp}
\end{equation}

\begin{theorem} Fix $\beta_k > 0$. If $g \in \Sigma_k(\beta_k)$ 
(of a regularity to be specified later), then there is a neighbourhood
$\calU$ of $g$ in the space of all conformally compact metrics 
on $M$ such that $\calU \cap \Sigma_k(\beta_k)$ is an analytic Banach 
submanifold of $\calU$, with respect to an appropriate Banach topology.
\end{theorem}

This theorem gives a rich class of conformally compact $\sk$-Yamabe
metrics on the manifold $M$, granting the existence of at least
one such metric. In particular, it states that the deformation
theory for this problem is always unobstructed whenever 
$\beta_k > 0$. As noted above, the analogue of this theorem
holds also when $M$ is compact without boundary, and the proof
is similar but even more straightforward. For the record, we
state this result too:
\begin{theorem} Fix $\beta_k > 0$. Let $g$ be a metric on the
compact manifold $M$ and $[g]$ its conformal class. Suppose that 
$\sk(A_g) = (-1)^k \beta_k$. Then there is a neighbourhood $\calU$ of 
$[g]$ in the space of conformal classes on $M$ such that every conformal 
class $[g']$ sufficiently near to $[g]$ contains a unique metric 
$g'_u = e^{2u}g'$ with $\sk(A_{g'_u}) = (-1)^k \beta_k$ which is
near to $g$; the set of these solution metrics fills out an (open piece
of an) analytic Banach submanifold, with respect to an appropriate Banach 
topology.
\end{theorem}

Let us return to conformally compact metrics, and connect Theorem~1 with 
the first theme discussed in the introduction. We recall that a 
Poincar\'e-Einstein metric is also a $\sk$-Yamabe metric for every $k = 1, 
\cdots,  n+1$. The constant $\beta_k$ must equal the constant
$\beta_k^0$ for hyperbolic space, 
so in particular the moduli space $\calE$ of P-E metrics is included in 
the intersection of the $\Sigma_k(\beta_k^0)$ over all $k$. However, 
more is true:

\begin{theorem} Within the class of conformally compact metrics
on $M$, there is an equality between the sets of Poincar\'e-Einstein
metrics and the metrics which are in $\Sigma_k(\beta_k^0)$ for every $k = 
1, \cdots, n+1$. In other words, with $\beta^0 = (\beta_1^0,
\cdots, \beta_{n+1}^0)$, 
\[
\Sigma(\beta^0) \equiv \bigcap_{k=1}^{n+1} \Sigma_k(\beta_k^0) = \calE.
\]
Hence $\calE$ is a finite intersection of locally closed Banach
submanifolds, and in particular is always closed in the space
of conformally compact metrics on $M$. 
\end{theorem}

Notice that if $\beta = (\beta_1, \cdots, \beta_{n+1})$
is any other $(n+1)$-tuple of numbers with $\beta_k > 0$, and
if $M$ is compact without boundary, then
\[
\Sigma(\beta) = \bigcap_{k=1}^{n+1} \Sigma_k(\beta_k)
\]
is typically empty. As we explain later, if $g$ is any metric in 
$\Sigma(\beta)$, then its Ricci tensor has constant eigenvalues.
In particular, metrics with $\nabla \Ric = 0$ are in $\Sigma(\beta)$
for some $\beta$. However, the reverse inclusion may not
be true and in any case is not well understood; only a few
partial results are known, e.g. in the K\"ahler case \cite{ADM}. 

The plan for the rest of this paper is as follows. \S 1 reviews
the structure of the functionals $\calF_k$ and their linearizations
$\calL_k$, and this is followed in \S 2 by a discussion of the function 
spaces and of the mapping properties of the $\calL_k$ on these
spaces. The deformation theory for the $\sk$-Yamabe equations
and the proof of Theorems~1 and 2 is the topic of \S 3. The (very
brief) proof of Theorem~3 and further discussion of some geometric
aspects of the moduli space $\calE$ is contained in \S 4.
Finally, \S 5 contains a list of some interesting open questions
raised by the results here.

We wish to thank Paul Yang for providing initial inspiration for these 
results during a fortuitous conversation at MSRI, during which he pointed 
out the dearth of examples of $\sigma_2$-Yamabe metrics in the negative case,
and suggested some possible lines of enquiry. Matt Gursky's encouragement 
was also quite helpful.

\section{The functionals $\calF_k$}

Let us fix a conformally compact metric $g_0$, which we may as
well take to be smooth, i.e.\  $g_0 = \rho^{-2}\olg_0$, where both $\rho$ 
and $\olg_0$ are ${\mathcal C}^\infty$ on $\Mbar$. Fix also a
constant $\beta_k > 0$. Recall that the metric $g = e^{2u}g_0$ is in 
$\Sigma_k(\beta_k)$, and so has $\sk(A_{g}) = (-1)^k \beta_k$, provided 
\[
\calF(g_0,u,\beta_k) = \sk(\nabla^2 u - du \otimes du + \frac12|\nabla u|^2 g_0
- A_{g_0}) - \beta_k e^{2ku} = 0.
\]
In this section we recall some facts about the ellipticity of this
operator and the structure of its linearization. These facts are
taken from \cite{V}, and we refer there for all proofs and further discussion.

To approach the issue of ellipticity, first consider the $k^{{\mathrm{th}}}$
elementary symmetric function $\sk$ as a function on vectors $\lambda = 
(\lambda_1, \cdots, \lambda_{n+1}) \in \RR^{n+1}$. Let $\Gamma_k^+$
denote the component of the open set $\{\lambda: \sk(\lambda) > 0\}$ 
containing the positive orthant $\{\lambda: \lambda_j > 0 \ \forall \, j\}$.
We note that these are all convex cones with vertices at the origin and
\[
\{\lambda: \lambda_j > 0\ \forall j\} = \Gamma_{n+1}^+
\subset \Gamma_n^+ \subset \cdots \subset \Gamma_1^+ = 
\{\lambda: \sigma_1(\lambda) > 0\}.
\]
Also, let $\Gamma_k^- = -\Gamma_k^+$. A real symmetric matrix $A$ is said to 
lie in $\Gamma_k^\pm$ if its eigenvalues lie in the corresponding set. 

\begin{proposition} 
If $A_g \in \Gamma_k^-$, then $u \to \calF_k(g,u,\beta_k)$ is elliptic at 
any solution of $\calF_k(g,u,\beta_k) = 0$.
\label{pr:ell}
\end{proposition}

The proof of this, in \cite{V}, relies on the computation of the linearization
of $\calF_k$ in the direction of the conformal factor $u$. The neatest 
formulation of this requires a definition from linear algebra.  For any 
real symmetric matrix $B$, and any $q = 0, \cdots, n+1$, define the 
$q^{{\mathrm{th}}}$ Newton transform of $B$ as the new (real, symmetric) 
matrix 
\[
T_q(B) = \sigma_q(B) I - \sigma_{q-1}(B) B + \cdots + (-1)^q B^q.
\]
Of course, $T_{n+1}(B) = 0$. Now suppose that $B = B(\e)$ depends 
smoothly on a parameter $\e$, and write $B'(0) = \dot{B}$. It is proved 
in \cite{R} that 
\[
\left. \frac{d\,}{d\e}\right|_{\e = 0} \sigma_k(B(\e)) = 
\langle T_{k-1}(B), \dot B \rangle_g.
\]

We apply this to the Schouten tensors associated to the family of metrics 
$g(\e) = e^{2\e \phi }g$ where $g \in \Sigma_k(\beta_k)$. We have
\[
B(\e) = -A_g + \e \nabla^2 \phi + \e^2 (\frac12 |\nabla \phi|^2 - 
d\phi \otimes d\phi)
\]
so that $B = -A_g$ and $\dot{B} = \nabla^2 \phi$. Hence
\begin{equation}
\left. D\calF_k\right|_{g,0}((0,\phi) \equiv
\calL_k \phi = \langle T_{k-1}(-A_g), \nabla^2 \phi \rangle - 2 k \beta_k \phi.
\label{eq:fdfk}
\end{equation}

The proof of Proposition~\ref{pr:ell} in general (i.e.\  when $g$ is
not necessarily a solution itself and when the linearization is computed
at some solution $u \neq 0$) relies on the
following ingredients: the convexity of $\Gamma_k^+$, the {\it concavity} 
of the function $\sk^{1/k}$ in $\Gamma_k^+$ and finally the
fact that $T_{k-1}(B)$ is positive definite when $B \in \Gamma_k^+$. 

Let us compute $\calL_k$ more explicitly when $g$ is hyperbolic,
or in fact, when $g$ is an arbitrary Poincar\'e-Einstein metric.
We shall always normalize the metric so the Einstein condition is
$\Ric = -n g$. Therefore, by (\ref{eq:AEg}), if $g$ is P-E then 
$A_g = -\frac12 g$, and so 
\[
T_{k-1}(-A_g) = 2^{1-k}T_{k-1}(g) = 2^{1-k}\sum_{j=0}^{k-1}(-1)^{j}
{n+1 \choose k-1-j} \equiv c_{k,n}.
\]
It is well-known, and easy to prove by induction, that
\[
c_{k,n} = 2^{1-k}{n \choose k};
\] 
in particular, this constant is always positive. Hence we obtain
the useful formula
\begin{equation}
\calL_k \phi = c_{k,n}\Delta \phi - 2k \beta_k \phi,
\label{eq:explk}
\end{equation}
which holds when $g$ is Poincar\'e-Einstein.

If $g \in \Sigma_k(\beta_k)$ is a more general solution (i.e.\  not
necessarily P-E), then $\calL_k$ is more complicated. However,
certain properties remain valid. 
\begin{proposition}
Suppose $g \in \Sigma_k(\beta_k)$, $\beta_k > 0$, and let $\calL_k$
denote the linearization of $\calF_k(g,u,\beta_k)$ at $u=0$.
Then 
\begin{equation}
\calL_k \phi = c_{k,n}\Delta_g \phi - 2k \beta_k \phi + \rho^3 E\phi,
\label{eq:formlk}
\end{equation}
where $E$ is a second order operator with bounded coefficients 
on $\Mbar$ (smooth if $\rho$ and $\olg = \rho^2 g$ are smooth),
and without constant term. Furthermore, if $\calL_k \phi = 0$
and $||\phi||_{L^\infty} < \infty$, then $\phi \equiv 0$.
\end{proposition}

The final statement follows directly from the asymptotic
maximum principle. A direct calculation yields the form of
$\calL_k$. 

We note that (\ref{eq:formlk}) may also be obtained from general principles 
involving the theory of uniformly degenerate operators \cite{Ma1}, 
\cite{Ma2}. Since some of the main results of this theory will be invoked 
later anyway, we digress briefly to explain this setup.
Choose coordinates $(x,y_1, \ldots, y_n)$, $x \geq 0$ near a point of the 
boundary of $\Mbar$. A second order operator $L$ is said to be uniformly 
degenerate if it may be expressed in the form
\begin{equation}
L = \sum_{j+|\alpha| \leq 2} a_{j,\alpha}(x,y)(x\del_x)^j (x\del_y)^\alpha.
\label{eq:unifdeg}
\end{equation}
The coefficients may be scalar or matrix-valued, and although we usually
assume they are smooth, it is easy to extend most of the main conclusions of 
this theory when they are polyhomogeneous, or of some finite regularity.
Operators of this type arise naturally in geometry, and in particular
all of the natural geometric operators associated to a conformally
compact metric are uniformly degenerate. Note that the error term 
$\rho^3 E$ in (\ref{eq:formlk}) is actually of the form $\rho E'$ where 
$E'$ is some second order uniformly degenerate operator without constant 
term.

The `uniformly degenerate symbol' of this operator is elliptic provided 
\[
\sigma(L)(x,y;\xi,\eta) = \sum_{j+|\alpha|=2} a_{j,\alpha}(x,y)
\xi^j \eta^\alpha \neq 0 \qquad \mbox{when}\quad (\xi,\eta) \neq 0. 
\]
(For systems, we require $\sigma(L)$ to be invertible as a matrix
when $(\xi,\eta) \neq 0$.) We also define the associated {\it normal operator} 
\[
N(L) = \sum_{j+|\alpha|\leq 2} a_{j,\alpha}(0,y)(s\del_s)^j (s\del_u)^\alpha.
\]
The boundary variable $y$ enters only as a parameter, while the `active' 
variables $(s,u)$ in this expression may be regarded as formal, but in
fact are naturally identified with linear coordinates on the inward pointing
half-tangent space $T^+_{(0,y)}M$. In particular
\begin{proposition}
If $g$ is a smooth conformally compact metric, then its Laplace-Beltrami
operator $\Delta_g$ is an elliptic uniformly degenerate operator
with normal operator
\begin{equation}
N(\Delta_g) = \Delta_{{\mathbb H}^{n+1}} = (s\del_s)^2 + s^2\Delta_u
-n s\del_s.
\label{eq:nolap}
\end{equation}
Furthermore, if $g \in \Sigma_k(\beta_k)$ for some $\beta_k > 0$,
then the linearization $\calL_k$ of $\calF_k(g,u,\beta_k)$ at 
$u=0$ is also elliptic and uniformly degenerate, with normal operator
\begin{equation}
N(\calL_k) = \calL_k^0 \equiv c_{k,n}((s\del_s)^2 + s^2\Delta_u
-n s\del_s) - 2k\beta_k.
\label{eq:nolinbk}
\end{equation}
\label{pr:strlk}
\end{proposition}

As we explain in the next section, the operator $\calL_k$ is 
is Fredholm on various natural function spaces. This specializes
a criterion which is applicable to other more general uniformly
degenerate operators $L$, namely that $L$ is Fredholm if and only 
if two separate ellipticity conditions hold: first, the symbol $\sigma(L)$
should be invertible, and in addition, the normal operator $N(L)$
must be invertible on certain weighted $L^2$ spaces.

\section{Function spaces and mapping properties}

Let $\calL_k$ be the linearization considered in the last section.
We shall now describe some of its mapping properties. As indicated
above, these properties also hold for more general elliptic, uniformly
degenerate operators $L$.

We first review one particular scale of function spaces which is convenient 
in the present setting, and then state the mapping properties on them
enjoyed by $\calL_k$. The material here is taken from \cite{Ma1}, to which
we refer for further discussion and proofs.

Fix a reference (smooth) conformally compact metric $g_0 = \rho^{-2}
\olg_0$; also, choose a smooth boundary coordinate chart $(x,y)$
as in the previous section, and recall the basic vector fields $x\del_x$ 
and $x\del_{y_j}$, $j = 1, \ldots, n$. Since $x$ is a smooth nonvanishing 
multiple of $\rho$ near $\dMbar$, these vector fields are all of uniformly 
bounded lengths with respect to $g_0$, and are also uniformly independent 
as $x \searrow 0$. There are two equivalent ways to define the
H\"older space $\Lambda^{0,\alpha}_0(M)$, $\ell \in {\mathbb N}$,
$0 < \alpha < 1$. In either case, it suffices to work in a boundary 
coordinate chart. The first is to set
\[
\Lambda^{0,\alpha}_0(M) = \left\{u: \sup \frac{|u(x,y) - u(x',y')|
(x+x')^\alpha}{|x-x'|^\alpha + |y-y'|^\alpha}\right\},
\]
where the supremum is taken first over all points $z=(x,y)$, $z'=(x',y')$,
$z \neq z'$, which lie in some coordinate cube $B$ centered at a point 
$z_0 = (x_0,y_0)$ of sidelength $\frac12 x_0$, and then over all such
cubes. The other is to let $B$ denote a ball of unit radius 
with respect to the metric $g_0$ centered at $z_0$, and to replace
the denominator in this definition by $\mbox{dist}_{g_0}\,(z,z')^\alpha$,
for $z,z' \in B$, and then take the same sequence of suprema.

This latter definition is more geometric, while the former definition
clearly implies the scale invariance of these spaces, namely that
if $u(z)$ is defined (and, say, compactly supported) in one of these 
coordinate charts and if we define $u_\e(z) = u(z/\e)$, then the associated 
norms of $u$ and $u_{\e}$ are the same. 

From here we can define a few other closely related spaces which
will be useful:
\begin{itemize}
\item For $\ell \in {\mathbb N}$, let 
\[
\Lambda^{\ell,\alpha}_0(M) = \left\{u: (x\del_x)^j (x\del_y)^\beta u 
\in \Lambda^{0,\alpha}_0(M) \ \forall \ j+|\beta| \leq \ell\right\}.
\]
\item For $0 \leq \ell' \leq \ell$, $\ell,\ell' \in \mathbb N$, let
\[
\Lambda^{\ell,\alpha,\ell'}_0(M) = \left\{u: \del_y^{\beta'}u \in
\Lambda^{\ell-|\beta'|,\alpha}_0(M)\ \mbox{for}\ |\beta'| \leq \ell'\right\}.
\]
\item Finally, for $\gamma \in \RR$, and $0 \leq \ell' \leq \ell$,
$\ell, \ell' \in {\mathbb N}$, let
\[
\rho^\gamma\Lambda^{\ell,\alpha,\ell'}_0(M) = \left\{u: 
u = \rho^\gamma \tilde{u},
\quad \mbox{where}\quad \tilde{u} \in \Lambda^{\ell,\alpha,\ell'}_0(M)\right\}.
\]
\end{itemize}

Thus the first of these are just the natural higher order H\"older spaces 
associated to the geometry of $g_0$, or equivalently, to differentiations
with respect to the vector fields $x\del_x$ and $x\del_y$. The second
of these allows up to $\ell'$ of the derivatives to be taken with respect
to the `nongeometric' vector fields $\del_y$, which are exponentially
large with respect to $g_0$ as $x \to 0$, but which are needed if one 
desires any sort of boundary regularity. The final spaces are just the 
usual weighted analogues. The corresponding norms are  
$||\cdot||_{\ell,\alpha}$, $||\cdot||_{\ell,\alpha,\ell'}$, and
$||\cdot||_{\ell,\alpha,\ell',\gamma}$, respectively.

We could equally easily have defined $L^2$-- and $L^p$--based Sobolev 
spaces, corresponding to differentiations with respect to the vector 
fields $x\del_x$ and $x\del_y$, as well as the corresponding partially 
tangentially regular and weighted versions. The mapping properties
we state below all have direct analogues for these spaces. However,
as usual the H\"older spaces are perhaps the simplest to deal
with for nonlinear PDE. 

As a further note, still fixing the reference metric $g_0$ we can use 
these definitions to define appropriate finite regularity spaces of 
tensor fields on $M$. In particular, we set, for $\gamma > 0$,
\[
{\mathfrak M}^{\ell,\alpha,\ell'}_\gamma(M) = \{g = g_0 + h: g > 0 
\ \mbox{and}\ h \in \rho^{-2+\gamma}\Lambda^{\ell,\alpha,\ell'}_0(M;S^2(M))\}.
\]

Now let us turn to the mapping properties of $\calL_k$. First of
all, it follows immediately from the definitions that 
\begin{equation}
\calL_k: \rho^\gamma \Lambda_0^{\ell+2,\alpha,\ell'}(M)
\longrightarrow \rho^\gamma\Lambda_0^{\ell,\alpha,\ell'}(M)
\label{eq:map}
\end{equation}
is a bounded mapping for any $\gamma \in \RR$ and $0 \leq \ell' 
\leq \ell$. However, this map is not well-behaved for many
values of the weight parameter $\gamma$.  There are two
ways this may occur. First if $\gamma$ is sufficiently large
positive, then it is not hard to see that (\ref{eq:map})
has an infinite dimensional cokernel, while dually, if $\gamma$
is sufficiently large negative, then (\ref{eq:map}) has
an infinite dimensional nullspace. Although we do not use
it here, less trivial is the fact that in either of these
two cases the mapping is semi-Fredholm (i.e.\ has closed
range and either the kernel or cokernel are finite-dimensional).

However, for certain values of $\gamma$ the range of this
mapping may not be closed. This is determined by a consideration
of the indicial roots of $\calL_k$. We say that $\gamma$
is an indicial root of $\calL_k$ if $\calL_k(\rho^\gamma) = 
O(\rho^\gamma+1)$ (note that because of the uniform degeneracy of 
$\calL_k$, $\calL_k(\rho^\gamma) = O(\rho^\gamma)$ is always
true). Thus $\gamma$ is an indicial root only if some special 
cancellation occurs. Using Proposition~\ref{pr:strlk}, it is
clear that the indicial roots of $\calL_k$ agree with those
of $\calL_k^0$, and then (\ref{eq:nolinbk}) shows that $\gamma$ 
is an indicial root if and only if
\[
c_{k,n}(\gamma^2 -n \gamma) - 2k\beta_k = 0,
\]
or in other words
\[
\gamma = \gamma_{\pm} = 
\frac{c_{k,n} n \pm \sqrt{(c_{k,n}n)^2 + 8k\beta_k c_{k,n}}}{2c_{k,n}}
= \frac{n}{2} \pm \sqrt{\frac{n^2}{4} + \frac{2k}{c_{k,n}}}.
\]
In particular
\[
\gamma_- < 0 < \gamma_+.
\]

The relevance of these indicial roots to the mapping properties
of (\ref{eq:map}) is that when $\gamma$ is equal to one of these
two values, then (\ref{eq:map}) does not have closed range.
At heart, this stems from the fact that the equation 
\[
\calL_k^0 u = s^{\gamma_\pm}
\]
has solution $u = c s^{\gamma_\pm}(\log s)$ for some constant $c$,
i.e.\ the inhomogeneous term is in the appropriate weighted H\"older
space but the solution $u$ just misses being in this space. 

Despite these cautions, we have the following basic result:

\smallskip

\noindent{\bf Mapping properties:} {\it If $\gamma_- < \gamma <
\gamma_+$, then the mapping (\ref{eq:map}) is Fredholm of index zero.}  

\smallskip

The main result of \cite{Ma1} is a considerably more general theorem 
of this sort for more general elliptic uniformly degenerate differential 
operators. There are two special features of $\calL_k$ which enter into 
the precise form of the statement here. First, there is a nontrivial
interval $(\gamma_- ,\gamma_+)$ between the two indicial roots 
$\gamma_\pm$, allowing for the possibility of a `Fredholm range'. 
Second, the Fredholm index is zero for $\gamma$ in this interval 
ultimately because $\calL_k$ is self-adjoint on $L^2(dV_g)$.

To show that $\calL_k$ is actually invertible when $\gamma$ is
in this Fredholm range, we note that at least when $\gamma \geq 0$
we could use the maximum principle to assert that the nullspace
is trivial, and then use the vanishing of the index to conclude
that the mapping is actually an isomorphism. In fact, this same 
reasoning may be used whenever $\gamma > \gamma_-$ on account
of another basic from \cite{Ma1}:

\smallskip

\noindent{\bf Regularity of solutions:} {\it If $\gamma_- <
\gamma < \gamma_+$ and $\phi \in x^{\gamma}\Lambda^{\ell+2,\alpha,
\ell'}_0$ is a solution of $\calL_k \phi = f$ where $f$ vanishes
to all orders at $\dMbar$, then as $x \to 0$, 
\[
\phi(x,y) \sim \sum_{j=0}^\infty \phi_j(y)x^{\gamma_+ + j},
\qquad \mbox{with} \qquad \phi_j(y) \in \calC^\infty(\del\Mbar),
\]
i.e.\  $\phi \in \rho^{\gamma_+}{\mathcal C}^\infty(\Mbar)$. 
}

\smallskip

In particular, if $\phi$ is in the nullspace of $\calL_k$ when 
$\gamma > \gamma_-$ then we may apply this (with $f=0$) to
obtain that $\phi = O(\rho^{\gamma_+})$ and so the maximum principle 
may be used as above to conclude that $\phi = 0$. 

\section{Perturbation theory in $\Sigma_k$}

Before preceding to the main deformation result, we review the
the structure of the space of (conformally compact) metrics on $M$.

There is a multiplicative action by scalar functions on
the set of conformally compact metrics:
\[
\rho^\gamma \Lambda^{\ell+2,\alpha,\ell'}_0(M) \times
{\mathfrak M}^{\ell+2,\alpha,\ell'}_\gamma(M) \longrightarrow 
{\mathfrak M}^{\ell+2,\alpha,\ell'}_\gamma(M),
\]
\[
(u,g) \longrightarrow e^{2u}g.
\]
The tangent space to this action at $(0,g)$ consists of the set 
\[
\left\{h \in \rho^\gamma\Lambda^{\ell+2,\alpha,\ell'}_0(M;S^2(T^*M)):
\tr_g(h) = 0\right\}.
\]
It is customary to model the set of conformal classes ${\mathfrak C}$ in a 
neighbourhood of $g$ by the set of tensor fields $h$ which are 
traceless, and also divergence free. The latter condition corresponds
to fixing a gauge transverse to the orbit of the action of the
diffeomorphism group; however, it is not really necessary for our
purposes to mod out by this action, and so we dispense with this gauge.

Now define the slice
\begin{equation}
{\mathfrak S}_\e = \left\{h \in \rho^\gamma\Lambda^{\ell+2,\alpha,\ell'}_0(
M;S^2(T^*M)):
\tr_g(h) = 0,\ ||h||_{\ell+2,\alpha,\ell',\gamma} < \e \right\},
\label{eq:slice}
\end{equation}
and identify a neighbourhood of $0$ in this space with the
space of conformal classes of conformally compact metrics
near $g$ (but not modulo diffeomorphisms!). 

\begin{theorem}
Fix $\beta_k > 0$ and suppose that $g \in \Sigma_k(\beta_k)$. Then there 
is a neighbourhood $\calU$ of $g$ in ${\mathfrak M}^{\ell+2,\alpha,\ell'}_\gamma$
such that $\calU \cap \Sigma_k(\beta_k)$ is a Banach manifold modelled on 
the slice ${\mathfrak S}_\e$ defined in (\ref{eq:slice}).
\end{theorem}

Having set things up carefully, the proof is almost immediate. Consider
the mapping
\begin{equation}
{\mathcal H}: 
{\mathfrak S}_\e \times \rho^\gamma\Lambda^{\ell+2,\alpha,\ell'}_0(M)
\longrightarrow {\mathfrak M}^{\ell,\alpha,\ell'}_\gamma
\label{eq:defH}
\end{equation}
defined by
\[
{\mathcal H}(h,u) = \calF_k(g+h,u,\beta_k).
\]
In a neighbourhood of $g$, the set $\Sigma_k(\beta_k)$ is identified with
the zero set of ${\mathcal H}$. In particular, $(h,u)=(0,0) \in
{\mathcal H}^{-1}(0)$.

To find all other nearby solutions, we shall apply the implicit function 
theorem, very much in the spirit of the closely related papers \cite{MS} 
and \cite{GL}. Thus we must check two things:
\begin{itemize}
\item The mapping ${\mathcal H}$ in (\ref{eq:defH}) is a smooth
mapping of open sets of Banach spaces, and
\item The linearization $\left. D {\mathcal H}\right|_{0,0}$ is
surjective between the appropriate tangent spaces.
\end{itemize}
The first of these is straightforward from the definitions and
(\ref{eq:trsch}), provided we choose the weight parameter $\gamma > 0$.
As for the other, recall that the restriction of this Frech\'et derivative 
to tangent vectors of the form $(0,\phi)$ corresponds to the 
operator $\calL_k$. We have already checked that this is
surjective on the tangent spaces provided we choose the
weight parameter $\gamma \in (\gamma_-,\gamma_+)$. Since $\gamma_+ > 0$,
these two restrictions on $\gamma$ are not inconsistent, and so we
now fix any $\gamma \in (0,\gamma_+)$. With this choice, we obtain
the existence of a smooth map 
\[
\Phi: {\mathfrak S}_\e 
\longrightarrow \rho^\gamma\Lambda^{\ell+2,\alpha,\ell'}_0(M)
\]
with $\Phi(0) = 0$, $\left. D\Phi\right|_0 = 0$ and such that
\[
{\mathcal H}(h,\Phi(h)) \equiv 0.
\]
Furthermore, all solutions of ${\mathcal H}(h,u)$ in a sufficiently
small neighbourhood of $(0,0$ are of this form. This concludes
the proof.  \hfill $\Box$

\smallskip

We omit the proof of Theorem 2 because it is nearly identical; 
indeed, the only difference is that standard elliptic theory
replaces the Fredholm theory for uniformly degenerate operators
we have quoted.

We have shown that $\Sigma_k(\beta_k)$ is a Banach submanifold in a 
neighbourhood of $g$, and furthermore that it may be regarded as a 
graph over the space of conformal classes, or at least those conformal 
classes near to $g$. For $k = 1$, every conformal class on $\Mbar$ 
contains a unique representative lying in $\Sigma_k(\beta_k)$, and
thus $\Sigma_1(\beta_1)$ is a graph globally over the space
$\mathfrak C$ of all conformal classes. It is not known whether
this remains true when $k > 1$, and thus we define
\begin{equation}
{\mathfrak C}_k = \left\{ {\mathfrak c} \in {\mathfrak C}: {\mathfrak c}
\ \mbox{contains at least one}\ g \in \Sigma_k(\beta_k)\right\}.
\label{eq:frakck}
\end{equation}
(Note we have suppressed the dependence of $\beta_k$ in this
notation because, if there is a solution $g \in {\mathfrak c}$ 
for the equation with one constant $\beta_k > 0$, then, simply
by scaling, we can produce a solution for every other positive constant.)
Note that ${\mathfrak C}_1 = \mathfrak C$, and Theorem~4 shows that 
${\mathfrak C}_k$ is open in ${\mathfrak C}$ for every $k$. 

\section{The moduli space $\calE$ of Poincar\'e-Einstein metrics}

We begin by recalling some previously known results about 
the moduli space ${\mathcal E}$ of Poincar\'e-Einstein metrics on
the manifold with boundary $M$, and then proceed to the
(nearly obvious) relationship of this space with the moduli spaces 
$\Sigma_k(\beta_k)$.

It is natural to formulate the existence theory of P-E metrics as an 
asymptotic boundary problem. Recall that any conformally compact 
manifold $g$ determines 
a conformal class $[\left.\rho^2 g\right|_{\dMbar}]$ on the
boundary of $\dMbar$. We then pose the following

\smallskip

\noindent{\bf Problem:} {\it Given a conformal class $[h_0]$ on
$\dMbar$, is it possible to find a Poincar\'e-Einstein
metric $g$ on $M$, such that the conformal class determined
by $g$ on $\dMbar$ is the specified class $[h_0]$?}

\smallskip

There are many particular solutions in the physics literature, cf.\ 
\cite{HP}, and this problem was considered on a formal level by 
Fefferman and Graham \cite{FG}. The first general existence result was 
obtained by Graham and Lee \cite{GL}, and consisted of a perturbation 
analysis about the standard hyperbolic metric on the ball $B^{n+1}$.
They showed that this problem is solveable for conformal classes 
sufficiently near the standard one. There are now several other 
related perturbation results \cite{Bi}, \cite{L}, and Anderson has 
announced a more global existence theory in four dimensions \cite{A3}. 

There are many interesting questions, including the description 
of asymptotic regularity of these metrics in terms of the regularity
of the conformal class $[h_0]$, problems involving uniqueness, etc.
Anderson \cite{A2} settles (appropriate versions of) these particular 
issues in the four dimensional case while \cite{GM} contains a different
approach to regularity and uniqueness valid in all dimensions.

The goal in this section is simply to provide a slightly different 
perspective onto this class of metrics. Fix the constants $\beta_k =
\beta_k^0$ corresponding to the standard hyperbolic metric,
and henceforth write $\Sigma_k(\beta_k^0)$ simply as $\Sigma_k$. 
In any case, for any $\ell$-tuple $J = \{j_1, \ldots, j_\ell\} \subset 
\{1, \ldots, n+1\}$, consider the intersection
\[
\Sigma_J = \Sigma_{j_1} \cap \cdots \cap \Sigma_{j_\ell},
\]
and 
\[ 
\Sigma = \Sigma_1 \cap \cdots \cap \Sigma_{n+1}.
\]

\begin{proposition} The moduli space of P-E metrics on $M$
agrees with this intersection of the submanifolds $\Sigma_k$:
\begin{equation}
{\mathcal E} = \Sigma. 
\label{eq:fint}
\end{equation}
\end{proposition}

The fact that ${\mathcal E} \subset \Sigma$ is obvious.
On the other hand, if $g \in \Sigma$, then we first notice that all 
the eigenvalues of $A_g$ (or equivalently,
of $\mbox{Ric}_g$) are constant; this is simply because they
are the roots of the characteristic polynomial, the coefficients
of which are precisely these higher ($\sk$) traces of $A_g$,
which are by assumption constant. To see that these constant
eigenvalues $\lambda_j$ agree with one another, and are equal to
the value on hyperbolic space, we can either proceed geometrically 
and note that near infinity (i.e.\  as $\rho \to 0$), the metric $g$ 
is asymptotically hyperbolic so we may evaluate the $\lambda_j$ in this 
limit, or else algebraically, by observing that by the particular 
choice of constants $\beta_k^0$, the characteristic
polynomial of $A_g$ must be simply that of hyperbolic space,
i.e.\  $(\lambda + 1/2)^{n+1}$. By either approach, we conclude
that $g$ is Einstein, and this gives the opposite inclusion.

There are various ways to interpret this equality of moduli spaces.

\begin{itemize}
\item Proposition~4 shows that the somewhat less tractable space 
${\mathcal E}$ is a finite intersection of submanifolds $\Sigma_j$,
each of which is an {\it analytic} submanifold, but more importantly, 
each of which has an unobstructed deformation theory. This amounts to 
some sort of figurative `factorization' of the Einstein equations into 
$n+1$ scalar (albeit fully nonlinear) equations.

\item Each $\Sigma_k$ is a (perhaps multiply sheeted) graph
over the open set of conformal classes ${\mathfrak C}_k$. But the
intersection $\Sigma$ lies over a much thinner, infinite codimension 
set ${\mathfrak C}_{\mathcal E} \subset {\mathfrak C}$. Elements of 
${\mathfrak C}_{\mathcal E}$ are precisely the conformal classes which 
contain P-E metrics! Another way to think of the basic asymptotic boundary 
problem for P-E metrics is that the set ${\mathfrak C}_{\mathcal E}(\dMbar)$ 
of conformal classes on $\dMbar$ for which this problem is solveable
coincides with the restrictions to $\dMbar$ of classes in  
${\mathfrak C}_{\mathcal E}$. This near-tautology presents an interesting 
(though probably intractable) way to think about the boundary problem:
one might first seek to extend the conformal class $[h_0]$ from the 
boundary to a `good' one (i.e.\  one in ${\mathfrak C}_{\mathcal E}$)
in the interior, and then consider the Einstein equation as 
a conformal problem within this class. 

\end{itemize}

\section{Open questions and further directions}

We conclude this note by raising a few other problems and
questions related to the results and methods here.

\begin{enumerate}

\item[a)] Because of the difficulty in obtaining a ${\mathcal C}^2$
estimate for the $\sk$-Yamabe problem when $k > 1$, it is worth wondering 
whether it might be worthwhile to pose a weaker version of this problem, 
at least for conformally compact metrics on manifolds with boundary: 
namely, given a conformal class $[h_0]$ on $\dMbar$, is it possible 
to extend this conformal class to at least some conformal class 
$[\olg]$ on the interior such 
that the $\sk$-Yamabe problem is solveable in $[\olg]$?
Probably there are infinitely many such extensions, as is
the case when $k=1$, but the added flexibility in this
formulation may be of some use. 

\item[b)] It seems central to understand whether ${\mathfrak C}_k =
{\mathfrak C}$, or in other words, whether every conformal class
on $M$ contains a conformally compact $\sk$-Yamabe metric.
Related to this is the observation that we do not know whether each 
of the submanifolds $\Sigma_k$ itself is closed; this depends 
ultimately on whether some version of this $\calC^2$ estimate holds.
It seems interesting, though, to ask which of the finite intersections
$\Sigma_J = \Sigma_{j_1} \cap \ldots \cap \Sigma_{j_\ell}$ are closed. 
Notice that if $1 \in J$, then this is
certainly true because the ${\mathcal C}^2$ estimate for the conformal 
factor is routine for the scalar curvature equation.

\item[c)] The regularity of the metrics $g \in \Sigma_k(\beta_k)$
is an interesting question. When $k=1$ this is resolved in 
\cite{Ma2}, cf. also \cite{MT}: if $g$ is a smooth conformally compact 
metric, then the conformal factor $u$ corresponding to the unique
solution $\tilde{g} = e^{2u}g \in [g]$ has a polyhomogeneous 
expansion. Presumably a similar result holds for all $k$. 
Note that unless $\gamma_+ \in {\mathbb N}$, this expansion
will involve nonintegral powers of $\rho$; this should not be viewed
negatively, since functions with expansions of this form
may be manipulated just as easily as smooth functions. 

\item[d)] The $\sk$-Yamabe problem considered here extends naturally
to the more general setting of the singular $\sk$-Yamabe problem:
given a smooth metric $g_0$ on a compact manifold $M$ and a closed subset 
$\Lambda \subset M$, when is it possible to find a conformally related 
metric $g = e^{2u}g_0$ which is both complete on $M \subset \Lambda$
and a $\sk$-Yamabe metric? When $k=1$ it is known that the dimension 
of $\Lambda$ is intimately related to the sign of the imposed scalar 
curvature of the solution, and very good existence results are known 
when $\Lambda$ is a submanifold \cite{AMc}, \cite{MP}. What is the
correct statement, and to what extent is this true when $k>1$?
There are a number of interesting analytic problems of this nature, 
and we shall return to this soon.

\item[e)] In general (not just in the conformally compact setting), 
$\calE$ sits inside the finite intersection $\cap \Sigma_k$.
Does it appear here as a finite codimensional analytic set, and if so,
is this related to some sort of Kuranishi reduction for the perturbation 
theory for $\calE$? 

\item[f)] Finally, it appears that very little is known about metrics 
with Ricci tensor having constant eigenvalues, but cf.\ \cite{ADM}. 
The (presumably) smaller class of metrics with parallel Ricci tensor is
more tractable, but it does not seem to be known if the eigenvalues
can be constant without the Ricci tensor being parallel. Examples would be 
very welcome. Also, in any setting (compact or conformally compact or ...)
it seems to be a very basic problem in Riemannian geometry to ask 
what are the possible $(n+1)$-tuples $(\beta_1,\ldots ,\beta_{n+1})$
for which $\Sigma_1(\beta_1) \cap \ldots \cap \Sigma_{n+1}(\beta_{n+1})$
is nonempty?

\end{enumerate}

\end{document}